\documentclass[12pt]{article}
\usepackage{amssymb,epsf}
\newtheorem{theorem}{Theorem}
\newtheorem{definition}{Definition}
\newtheorem{statement}{Statement}
\newtheorem{lemma}{Lemma}
\newtheorem{corollary}{Corollary}
\newcommand{\kC}{{\mathbb C}}
\newcommand{\kQ}{{\mathbb Q}}
\newcommand{\proj}[1]{{\mathbb P}^{#1}}
\newcommand{\begproof}{\ensuremath{\blacktriangleleft}}
\newcommand{\enproof}{\ensuremath{\blacktriangleright}}
\usepackage{maplestd2e}

\begin{document}
{\center \LARGE Calculation of two Belyi pairs.\\[3mm]
\normalsize Dremov Vladimir //21.11.2008\\[1cm]
}

\section{Introduction}
Our calculation is a part of the work which results in a catalogue \cite{maxplank}.

We want to find the Belyi pair of the following dessins (opposite sides are identified
to obtain an oriented surface of genus one):
 \epsfysize=1cm\epsfbox{pict.1}  \epsfysize=1cm\epsfbox{pict.2}

Our calculation uses properties of Mulase-Penkava differential (see below).
Instead of solving equations for Belyi pair, we construct equations for two quadratic
differentials constructed from the Belyi function.

We construct appropriate coordinates. Then we solve equations to find coefficients
of a curve equation and of the polynomials corresponding to the considered
quadratic differentials.

\section{Notations}
We use valencies of the dessin and write the structure of divisors $(\beta)$
and $(\beta-1)$. We find all of the dessins with the same valencies. There
are two dessins with these valencies. We obtain the following result:
\begin{statement}
\label{defdess}
There exists two non-isomorphic Belyi pairs $(X,\beta)$ of genus $1$ with divisorial structure
$$(\beta)=5A_1+3A_2-7C_1-C_2\eqno$$ and $$(\beta-1)=2(B_1+B_2+B_3+B_4)-7C_1-C_2.$$ There are
the only two isomorphism classes of such dessins, namely,
$$\epsfysize=3cm\epsfbox{pict.1}\epsfysize=3cm\epsfbox{pict.2}$$
\end{statement}

\begin{definition}
The Belyj function (and the corresponding dessin) are called \emph{clean} if
all preimages in $\beta^{-1}(1)$ are double points \cite{LSchneps}.
\end{definition}

\begin{definition}
Mulase-Penkava operator $MP$ is a differential operator, associating the meromorphic
quadratic differential to a nonconstant meromorphic function on a fixed Riemann surface
$$MP:f\mapsto \frac{1}{4\pi^2}\frac{(df)^2}{f(1-f)}$$
\end{definition}
This operator is considered in the article \cite{Mul_Penk}. We use the following elementary properties of $MP(f)$:

\begin{lemma}
If $f$ is a clean Belyi function and $(f)=\sum_{i=1}^p a_iA_i-\sum_{j=1}^q c_jC_j$ then
$(MP(f))=\sum_{i=1}^p (a_i-2)A_i-\sum_{j=1}^q 2C_j$
\end{lemma}
\begin{lemma}
If $f$ has a pole of order $k$ at the point $C$ then the residue of $MP(f)$ at the point $C$ is equal to
$$Res_C(MP(f))=-\frac{k^2}{4\pi^2}.$$
\end{lemma}

We calculate the two Belyi pairs from the statement \ref{defdess} above (the result stated in the theorem \ref{BPairs},
p.~\pageref{BPairs}). We use $MP(\beta)$ and $MP(\beta^{-1})$ in our calculations.

Denote $(X,\beta)$ one of the two Belyi pairs, satisfying the above divisorial relations.

The dessin $\beta^{-1}[0,1]$ is clean. The points $A_1,A_2,C_1,C_2 \in X$ are
uniquely defined for this dessin by $(\beta)=5A_1+3A_2-7C_1-C_2$.

\section{Definition of coordinates}

There is the only elliptic involution $\tau:X\to X$ for which $\tau(C_1)=C_2$.

Consider coordinates $x:X\to X/_{<1,\tau>}\simeq \proj{1}(\kC)=\kC \cup \{\infty\}$ for which
$$x(C_1)=x(C_2)=\infty,$$
$$x(A_1)=0.$$
A projection $x$ of degree $2$ is defined up to multiplication $x\mapsto kx$, $k\in \kC$, $k\not=0$. Denote
$W_1,W_2,W_3,W_4$ the critical points of $x$.

Consider coordinates $y:X\to \proj{1}(\kC)=\kC \cup \{\infty\}$ for which
$(y)=W_1+W_2+W_3+W_4-2C_1-2C_2$. A coordinate $y$ has degree $4$ and is defined
up to multiplication $y\mapsto ky$, $k\in \kC$, $k\not=0$.

\begin{lemma}$y(A_1)\not=0$\end{lemma}
\begproof One easily checks that $(MP(\beta))=3A_1+A_2-2C_1-2C_2\equiv 0$.

Suppose that $y(A_1)=0$.\\
Hence $A_1\in \{W_1,W_2,W_3,W_4\}$ and $2A_1\equiv C_1+C_2$. \\
Hence $3A_1+A_2-2C_1-2C_2\equiv 3A_1+A_2-2\cdot(C_1+C_2) \equiv 3A_1+A_2-4A_1\equiv A_2-A_1 \equiv 0$.
Since $A_1\not=A_2$, we have come to a contradiction.
\enproof

\vspace{1cm}
Let us normalize the coordinate $y$ by the condition $y(A_1)=1$.

Let us consider the coordinates $x$ for which $\frac{x^2}{y}(C_1)=1$ (there are two such coordinates).
\section{Working with equations}
\subsection{Equation of the curve}
We know $x(A_1)=0, y(A_1)=1$ and $\frac{x^2}{y}(C_1)=1$. Hence the curve $X$ can be defined
by an equation of the form $$y^2=f(x):=1+ax+bx^2+cx^3+x^4.$$
 The coordinate $x$ is defined up to a sign.
\subsection{Equations for $MP(\beta)$}
Consider ${\omega}=\frac{dx}{y}=2\frac{dy}{f'}$.

Define the function $u$ by the equation $MP(\beta)=:u\cdot {\omega}^2$.

$$(u)=3A_1+A_2-2C_1-2C_2.$$

The function $u$ can be written in a following form:
$$u=p+qx+rx^{2}+sy$$

\subsubsection{Equations at $A_1$}
At the point $A_1$ (recall $x(A_1)=0$, $y(A_1)=1$) we have
$$y=\,1+\frac{1}{2}\,ax+ \left( \frac{1}{2}\,b-\frac{1}{8}\,{a}^{2} \right) {x}^{2}+O(x^3);$$
Substituting in $u=O(x^3)$, we obtain
$$\left( r+\frac{1}{2}\,sb-\frac{1}{8}\,s{a}^{2} \right) {x}^{2}+
\left( q+\frac{1}{2}\,sa \right) x+p+s+O(x^3)=O(x^3).$$

We solve this and obtain the following:
$$\left\{
p=-s,q=-\frac{1}{2}\,sa,r=-\frac{1}{2}\,sb+\frac{1}{8}\,s{a}^{2}\right\}$$

\subsubsection{Equations for quadratic residues of $MP(\beta)$}
$MP(\beta)=(r+s+o(1))\frac{dx^2}{x^2}$ at the point $C_1$;\\
$MP(\beta)=(r-s+o(1))\frac{dx^2}{x^2}$ at the point $C_2$. Hence

$$\frac{Res_{C_1}[MP(\beta)]}{Res_{C_2}[MP(\beta)]}=\frac{r+s}{r-s}=49;$$
$$\left\{ \frac{p}s=-1,\frac{q}s=-\frac{1}{2}\,a,\frac{r}s={\frac {25}{24}},b=\frac{1}{4}\,{a}^{2}-{\frac {25}{12}} \right\}.$$

\[f\, = \,1+ax+ \left( \frac14\,{a}^{2}-{\frac {25}{12}} \right) {x}^{2}\\
\mbox{}+c{x}^{3}+{x}^{4}\]

\subsubsection{Coordinates of $A_2$}
Coordinates of the point $A_2$ form a solution of the equations $y^2=f(x)$ and $u(x,y)=0$.
Solving this system we obtain:

$$A_2=(x={\frac {576}{49}}\,c+{\frac {600}{49}}\,a,y=1-{\frac {705888}{2401}}\,ac
-{\frac {360300}{2401}}\,{a}^{2}-{\frac{345600}{2401}}\,{c}^{2})$$

\begin{lemma}$x(A_2)\not=0$ and $25a+24c\not=0$.
\end{lemma}
\begproof One could easily check the following sequence of statements:\\
$x(A_2)=0\Leftrightarrow 25a+24c=0 \Rightarrow (x(A_2)=0,y(A_2)=1) \Rightarrow A_2=A_1$.

Final equation is false. Hence $x(A_2)=0$ and $25a+25c=0$ are also impossible.
\enproof

\subsection{Equations for $MP(\beta^{-1})$}
\begin{lemma}
$MP(\beta^{-1})=\frac{MP(\beta)}{(-\beta)}$
\end{lemma}
\begproof Both sides of the equation are equal to
$\frac{1}{4\pi^2}\frac{(d\beta)^2}{\beta^2(\beta-1)}$.
\enproof

\begin{corollary}
$$\left(MP(\beta^{-1})\right)=-2A_1-2A_2+5C_1-C_2.$$
\end{corollary}

Recall ${\omega}=\frac{dx}{y}=2\frac{dy}{f'}$.

\begin{lemma}\label{QRNot}$MP(\beta^{-1})$ can be written in the following form:
$$MP(\beta^{-1})=\frac{{\omega}^2}{Q(x)+yR(x)}$$ for
$Q(x)=q_{{0}}+q_{{1}}x+q_{{2}}{x}^{2}+q_{{3}}{x}^{3}+q_{{4}}{x}^{4}+q_{{5}}{x}^{5}$,\\
$R(x)=r_{{0}}+r_{{1}}x+r_{{2}}{x}^{2}+r_{{3}}{x}^{3}$, $q_j,r_j \in \kC$.
\end{lemma}
\begproof
Function $\frac{{\omega}^2}{MP(\beta^{-1})}$ has a divisor
$(\frac{{\omega}^2}{MP(\beta^{-1})})=2A_1+2A_2-5C_1+C_2$.
Hence $\frac{{\omega}^2}{MP(\beta^{-1})}\in L(5C_1+5C_2)$
and this function can be written as a polynomial in $\kC(x,y)$
of the form $Q(x)+yR(x)$ as stated above.
\enproof

\begin{definition} $Q\in \kC[x]$, $R\in \kC[x]$,
$q_j\in \kC$, $r_j\in \kC$ are the polynomials
and their coefficients introduced in the lemma {\rm \ref{QRNot}}.
\end{definition}

\subsubsection{Equations at $C_2$}
At the point $C_2$ we have for $y=y(x)$ the following approximation:
$yx^3+O(x^{-1})=-{x}^{5}-\frac12\,c{x}^{4}+{\frac {25}{24}}\,{x}^{3}-\frac18\,{a}^{2}{x}^{3}+\frac18\,{c}^{2}{x}^{3}-
\frac{1}2\,a{x}^{2}-{\frac {25}{48}}\,c{x}^{2}
+\frac1{16}\,c{a}^{2}{x}^{2}-\frac1{16}\,{c}^{3}{x}^{2}+{\frac {1}{128}}\,{a}^{4}x+{\frac {25}{64}}\,{c}^{2}x
+{\frac {5}{128}}\,{c}^{4}x-{\frac {25}{192}}\,{a}^{2}x+\frac{1}4\,acx+{\frac {49}{1152}}\,x
-{\frac {3}{64}}\,{c}^{2}{a}^{2}x+{\frac {5}{128}}\,{c}^{3}{a}^{2}-{\frac {3}{256}}\,c{a}^{4}
+\frac1{16}\,{a}^{3}-{\frac {25}{48}}\,a-{\frac {433}{768}}\,c-\frac3{16}\,{c}^{2}a
-{\frac {125}{384}}\,{c}^{3}-{\frac{7}{256}}\,{c}^{5}
+{\frac {25}{128}}\,c{a}^{2}+O(x^{-1})$.

Substituting in $Q(x)+yR(x)=O(x^{-1})$ (see lemma \ref{QRNot}, p.~\pageref{QRNot}), we obtain:\\
$
\, \left( q_{{5}}-r_{{3}} \right) {x}^{5}+ \left( -r_{{2}}
-\frac{1}{2} cr_{{3}}+q_{{4}} \right) {x}^{4}+$\\
$
+ \left( -\frac{1}{8} {a}^{2}r_{{3}}-\frac{1}{2} cr_{{2}}+{\frac {25}{24}}\,r_{{3}}
-r_{{1}}+q_{{3}}+\frac{1}{8} {c}^{2}r_{{3}}\right) {x}^{3}+$\\
$
+( -\frac{1}{2} cr_{{1}}-{\frac {25}{48}}\,cr_{{3}}
+\frac{1}{8} {c}^{2}r_{{2}}+q_{{2}}+{\frac {25}{24}}\,r_{{2}}
-\frac{1}{8} {a}^{2}r_{{2}}-r_{{0}}+\frac{1}{16}\,c{a}^{2}r_{{3}}
-\frac{1}{2} ar_{{3}}-\frac{1}{16}\,{c}^{3}r_{{3}} ) {x}^{2}+$\\
$
+ ( \frac{1}{8} {c}^{2}r_{{1}}+\frac{1}{16}\,c{a}^{2}r_{{2}}+q_{{1}}+{\frac {25}{24}}\,r_{{1}}-{\frac {25}{48}}\,cr_{{2}}
+{\frac {1}{128}}\,{a}^{4}r_{{3}}-{\frac {3}{64}}\,{c}^{2}{a}^{2}r_{{3}}-\frac{1}{2} cr_{{0}}
-\frac{1}{16}\,{c}^{3}r_{{2}}-{\frac {25}{192}}\,{a}^{2}r_{{3}}+{\frac {25}{64}}\,{c}^{2}r_{{3}}
-\frac{1}{8} {a}^{2}r_{{1}}-\frac{1}{2} ar_{{2}}+{\frac {49}{1152}}\,r_{{3}}+
\frac{1}{4} acr_{{3}}+{\frac {5}{128}}\,{c}^{4}r_{{3}}) x+$\\
$
+{\frac {5}{128}}\,{c}^{3}{a}^{2}r_{{3}}-{\frac {3}{256}}\,c{a}^{4}r_{{3}}
+\frac{1}{4} acr_{{2}}-{\frac {3}{64}}\,{c}^{2}{a}^{2}r_{{2}}+{\frac {25}{128}}\,c{a}^{2}r_{{3}}
-{\frac {25}{48}}\,cr_{{1}}+\frac{1}{16}\,c{a}^{2}r_{{1}}-\frac{3}{16}\,{c}^{2}ar_{{3}}+q_{{0}}+{\frac {25}{24}}\,r_{{0}}
+{\frac {49}{1152}}\,r_{{2}}+{\frac {5}{128}}\,{c}^{4}r_{{2}}-{\frac {25}{192}}\,{a}^{2}r_{{2}}
+\frac{1}{16}\,{a}^{3}r_{{3}}-{\frac {125}{384}}\,{c}^{3}r_{{3}}-{\frac {7}{256}}\,{c}^{5}r_{{3}}
-\frac{1}{8} {a}^{2}r_{{0}}+\frac{1}{8} {c}^{2}r_{{0}}-{\frac {433}{768}}\,cr_{{3}}+{\frac {25}{64}}\,{c}^{2}r_{{2}}
-\frac{1}{16}\,{c}^{3}r_{{1}}+{\frac {1}{128}}\,{a}^{4}r_{{2}}-{\frac {25}{48}}\,ar_{{3}}
-\frac{1}{2} ar_{{1}}+O(x^{-1})=O(x^{-1})$.

We solve this and obtain the formulas for $q_0,q_1,q_2,q_3,q_4,q_5$ as polynomials in $\kQ[a,c,r_0,r_1,r_2,r_3]$:\\
$q_{{5}}=r_{{3}}$,\\
$q_{{4}}=r_{{2}}+\frac{1}{2} cr_{{3}}$,\\
$q_{{3}}=\frac{1}{8} {a}^{2}r_{{3}}+\frac{1}{2} cr_{{2}}-{\frac {25}{24}}\,r_{{3}}
+r_{{1}}-\frac{1}{8} {c}^{2}r_{{3}}$,\\
$q_{{2}}=\frac{1}{2} cr_{{1}}+{\frac {25}{48}}\,cr_{{3}}-\frac{1}{8} {c}^{2}r_{{2}}
-{\frac {25}{24}}\,r_{{2}}+\frac{1}{8} {a}^{2}r_{{2}}+r_{{0}}-\frac{1}{16} c{a}^{2}r_{{3}}
+\frac{1}{2} ar_{{3}}+\frac{1}{16} {c}^{3}r_{{3}}$,\\
$q_{{1}}=-\frac{1}{8} {c}^{2}r_{{1}}-\frac{1}{16} c{a}^{2}r_{{2}}-{\frac {25}{24}}\,r_{{1}}
+{\frac {25}{48}}\,cr_{{2}}-{\frac {1}{128}}\,{a}^{4}r_{{3}}+{\frac {3}{64}}\,{c}^{2}{a}^{2}r_{{3}}
+\frac{1}{2} cr_{{0}}+\frac{1}{16} {c}^{3}r_{{2}}+{\frac {25}{192}}\,{a}^{2}r_{{3}}-{\frac {25}{64}}\,{c}^{2}r_{{3}}
+\frac{1}{8} {a}^{2}r_{{1}}+\frac{1}{2} ar_{{2}}-{\frac {49}{1152}}\,r_{{3}}-\frac{1}{4} acr_{{3}}
-{\frac {5}{128}}\,{c}^{4}r_{{3}}$,\\
$q_{{0}}=-\frac{1}{16} {a}^{3}r_{{3}}+{\frac {125}{384}}\,{c}^{3}r_{{3}}
+{\frac {7}{256}}\,{c}^{5}r_{{3}}-{\frac {1}{128}}\,{a}^{4}r_{{2}}+\frac{1}{8} {a}^{2}r_{{0}}
-{\frac {25}{24}}\,r_{{0}}-{\frac {49}{1152}}\,r_{{2}}+{\frac {3}{64}}\,{c}^{2}{a}^{2}r_{{2}}
+\frac{1}{2} ar_{{1}}-{\frac {5}{128}}\,{c}^{4}r_{{2}}+{\frac {3}{256}}\,c{a}^{4}r_{{3}}
+{\frac {25}{192}}\,{a}^{2}r_{{2}}+\frac{3}{16} {c}^{2}ar_{{3}}-{\frac {25}{128}}\,c{a}^{2}r_{{3}}
-\frac{1}{4} acr_{{2}}+{\frac {25}{48}}\,cr_{{1}}-\frac{1}{16} c{a}^{2}r_{{1}}
+\frac{1}{16} {c}^{3}r_{{1}}-{\frac {5}{128}}\,{c}^{3}{a}^{2}r_{{3}}
+{\frac {433}{768}}\,cr_{{3}}+{\frac {25}{48}}\,ar_{{3}}-\frac{1}{8} {c}^{2}r_{{0}}
-{\frac {25}{64}}\,{c}^{2}r_{{2}}$.

\subsubsection{Equations at $A_1$}
At the point $A_1$ we have $y=1+\frac12 ax+O(x^2)$.

Substituting in $Q(x)+yR(x)=O(x^2)$ (see lemma \ref{QRNot}, p.~\pageref{QRNot}), we obtain the following:

$-96\,r_{{1}}x-96\,r_{{0}}-98\,r_{{2}}+1152\,axr_{{0}}+1200\,xcr_{{2}}-
98\,xr_{{3}}-144\,xc{a}^{2}r_{{2}}-576\,acr_{{2}}-576\,xacr_{{3}}+108
\,{c}^{2}{a}^{2}r_{{2}}+108\,x{c}^{2}{a}^{2}r_{{3}}-90\,{c}^{3}{a}^{2}
r_{{3}}+27\,c{a}^{4}r_{{3}}-144\,c{a}^{2}r_{{1}}+432\,{c}^{2}ar_{{3}}-
450\,c{a}^{2}r_{{3}}-288\,x{c}^{2}r_{{1}}+1200\,cr_{{1}}+1152\,xar_{{2
}}-900\,x{c}^{2}r_{{3}}-900\,{c}^{2}r_{{2}}+144\,{c}^{3}r_{{1}}+288\,x
{a}^{2}r_{{1}}-18\,{a}^{4}r_{{2}}+144\,x{c}^{3}r_{{2}}+1200\,ar_{{3}}-
18\,x{a}^{4}r_{{3}}+1152\,ar_{{1}}-90\,{c}^{4}r_{{2}}-90\,x{c}^{4}r_{{
3}}+300\,{a}^{2}r_{{2}}+300\,x{a}^{2}r_{{3}}-144\,{a}^{3}r_{{3}}+750\,
{c}^{3}r_{{3}}+63\,{c}^{5}r_{{3}}+1152\,xcr_{{0}}+288\,{a}^{2}r_{{0}}-
288\,{c}^{2}r_{{0}}+1299\,cr_{{3}}+O(x^2)=O(x^2)$.

\subsubsection{Equations at $A_2$}
We assume further that $y(A_2)\not=0$. We do not include the details
of calculations for the case $y(A_2)=0$ as they give no solutions.

At the point $A_2$ we have\\
$y= \left( 1-{\frac {705888}{2401}}\,ac-{\frac {360300}{2401}}\,{a}^{2}-{
\frac {345600}{2401}}\,{c}^{2} \right)\times$\\
$\times ( 1+{\frac {147}{2}}\,
\frac
{( -1920800\,c-1961617\,a+847310496\,c{a}^{2}+288240100\,{a}^{3}+271060992
\,{c}^{3}+830131200\,{c}^{2}a )}
{( -2401+705888\,ac+360300\,{a}^{2}+345600\,{c}^{2} ) ^{2}}
( x-{\frac {576}{49}}\,c-{\frac {600}{49}}\,a)
 )
+O\left(\left( x-{\frac {576}{49}}\,c-{\frac {600}{
49}}\,a \right) ^{2}\right) $

As in the previous case, substituting it into
$$Q(x)+yR(x)=O\left(\left( x-{\frac {576}{49}}\,c-{\frac {600}{49}}\,a \right) ^{2}\right),$$
 we obtain the equations.

\subsubsection{Solving equations at $A_1$ and at $A_2$}
For variables $r_0,r_1,r_2,r_3$ we have obtained
the $2$ linear homogeneous equations at $A_1$ and the 2 linear homogeneous equations at $A_2$.
The solution $r_j=0$ is inconsistent with our assumptions. Hence the determinant
of this system is equal to zero. The determinant can be written as follows:\\
\label{det_equation}$\,29365647704064\, ( -24786\,{a}^{9}c+2594160450\,{a}^{2}-103766418\,ac-2490394032\,{c}^{2}+60289110\,{c}^{4}
-367482654\,{a}^{4}+323616384\,{c}^{2}{a}^{2}-22842\,{c}^{9}a+1913139\,{c}^{8}-793881\,{a}^{8}
+12150\,{a}^{10}+11664\,{c}^{10}+113888592\,{c}^{3}{a}^{3}+28335096\,{c}^{4}{a}^{2}-20615148\,c{a}^{5}
-84035232\,{c}^{2}{a}^{4}+419560344\,{a}^{3}c-435983184\,a{c}^{3}-95264100\,a{c}^{5}+22690800\,{a}^{6}
+34999992\,{c}^{6}-5596290\,{a}^{4}{c}^{4}-2150064\,{a}^{5}{c}^{3}+6627096\,{a}^{3}{c}^{5}+1959552\,{a}^{6}{c}^{2}
+2517480\,{a}^{2}{c}^{6}-5692032\,a{c}^{7}+1215000\,{a}^{7}c-142884\,{c}^{5}{a}^{5}+20412\,{c}^{4}{a}^{6}
+27216\,{c}^{6}{a}^{4}+97200\,{c}^{3}{a}^{7}+93312\,{c}^{7}{a}^{3}-36450\,{c}^{8}{a}^{2}-34992\,{c}^{2}{a}^{8}-13841287201
)\times$\\
$\times \left(-2401+705888\,ac+360300\,{a}^{2}+345600\,{c}^{2} \right)\times$\\
$\times  \left(24\,c+25\,a \right) ^{4}$.

\subsection{Equations for the quadratic residues of $MP(\beta^{-1})$}
$MP(\beta^{-1})=\frac{{\omega}^2}{Q(x)+yR(x)}=\frac{dx^2}{y^2(Q(x)+yR(x))}$

\label{kdefsection}

Denote $k_1:=\frac{Q(x)+yR(x)}{x^2}(A_1)$ and $k_2:=\frac{Q(x)+yR(x)}{(x-x(A_2))^2}(A_2)$.

Using these notations, we can write an equation as follows:
$$\left(\frac{y(A_2)}{y(A_1)}\right)^2 \frac{k_2}{k_1}=\frac{25}{9}.$$

\subsubsection{Finding coefficient $k_1$ at $A_1$}
We use an approximation $y=\,1+\frac12 ax-{\frac {25}{24}}\,{x}^{2}+O(x^3)$ and
find $k_1$ as a coefficient of $Q(x)+yR(x)$ at $x^2$

$k_1=\frac {1}{2352}\, ( \left( 24\,c+25\,a \right)
( 243\,{c}^{8}-648\,{a}^{2}{c}^{6}+5400\,{c}^{6}+7776\,a{c}^{5}
+486\,{a}^{4}{c}^{4}+64854\,{c}^{4}-8100\,{c}^{4}{a}^{2}+129600\,a{c}^{3}
-15552\,{c}^{3}{a}^{3}+345600\,{c}^{2}+705888\,ac-129600\,{a}^{3}c
+7776\,c{a}^{5}-81\,{a}^{8}+2700\,{a}^{6}
+360300\,{a}^{2}-64854\,{a}^{4}-2401 )
)(
27\,{c}^{6}+1755\,{c}^{4}-81\,{c}^{4}{a}^{2}-864\,a{c}^{3}+14697\,{c}^{2}+81\,{c}^{2}{a}^{4}-2646\,{c}^{2}{a}^{2}+864\,{a}^{3}c-288\,ac\\
\mbox{}-14409\,{a}^{2}+891\,{a}^{4}+2401-27\,{a}^{6})^{-1}$

\subsubsection{Finding coefficient $k_2$ at $A_2$}
We use a Taylor approximation for $y$ at $A_2$ and
find $k_2$ as a coefficient of $Q(x)+yR(x)$ at $\left(x-x(A_2)\right)^2$
(much like to the previous section, but longer formulas).
\subsubsection{The residue equation}
The equation considered is $9k_2 y(A_2)^2-25k_1=0$ (see section \ref{kdefsection},
p.~\pageref{kdefsection}), the left side is a rational function in $\kQ(a,c)$.
Numerator of this function is a polynomial of degree $21$ in each of variables $a,c$
and looks like follows (all degrees of monomials in the second factor are of even degree):\\
$\left( 25\,a+24\,c \right) (76527504000000\,{a}^{20}-468348324480000\,{a}^{19}c$\\
$+676870467379200\,{a}^{18}{c}^{2}+1542912026886144\,{a}^{17}{c}^{3}$\\
\dots\\
$-628201913088647840269231422\,a{c}^{3}-449810685236937774900955707\,{c}^{4}$\\
$-876053650539179213151415953\,{a}^{2}-1692722483624861493698125134\,ac$\\
$-812124909439798006329946263\,{c}^{2}+7819771121260579336605617)$.

\subsection{Solving equations}
The determinant on the page \pageref{det_equation} has the following $3$ factors:
\begin{enumerate}
\item $29365647704064 \left(24\,c+25\,a \right) ^{4}$
\item $\left(-2401+705888\,ac+360300\,{a}^{2}+345600\,{c}^{2} \right)$
\item $( -24786\,{a}^{9}c+2594160450\,{a}^{2}-103766418\,ac-2490394032\,{c}^{2}+60289110\,{c}^{4}
-367482654\,{a}^{4}+323616384\,{c}^{2}{a}^{2}-22842\,{c}^{9}a+1913139\,{c}^{8}-793881\,{a}^{8}
+12150\,{a}^{10}+11664\,{c}^{10}+113888592\,{c}^{3}{a}^{3}+28335096\,{c}^{4}{a}^{2}-20615148\,c{a}^{5}
-84035232\,{c}^{2}{a}^{4}+419560344\,{a}^{3}c-435983184\,a{c}^{3}-95264100\,a{c}^{5}+22690800\,{a}^{6}
+34999992\,{c}^{6}-5596290\,{a}^{4}{c}^{4}-2150064\,{a}^{5}{c}^{3}+6627096\,{a}^{3}{c}^{5}+1959552\,{a}^{6}{c}^{2}
+2517480\,{a}^{2}{c}^{6}-5692032\,a{c}^{7}+1215000\,{a}^{7}c-142884\,{c}^{5}{a}^{5}+20412\,{c}^{4}{a}^{6}
+27216\,{c}^{6}{a}^{4}+97200\,{c}^{3}{a}^{7}+93312\,{c}^{7}{a}^{3}-36450\,{c}^{8}{a}^{2}-34992\,{c}^{2}{a}^{8}
-13841287201)$.
\end{enumerate}

The first factor is nonzero.

The system of the second factor and the residue equation is
inconsistent with our assumptions ($x(A_2)\not=0$ and $y(A_2)\not=0$).

The system of the third factor and the residue equation gives us the desired solutions as shown below.
Resultant of those $2$ equations is equal to zero and has the following factors:
\begin{enumerate}
\item $301 a^4+2688 a^2+36864$
\item $133 a^4+896 a^2-12288$
\item $4725 a^4+342405 a^2+4477456$
\item $5670 a^4-1439865 a^2+13942756$
\item $190005517894500 a^8+36552364751718900 a^6+7708662622309824945 a^4
+69471491411890643040 a^2+1517090351363521026304$
\end{enumerate}
Hence one of the factors 1.-5. is equal to zero. We consider the corresponding cases.

The first and the second cases are inconsistent with assumptions $x(A_2)\not=0$ and $y(A_2)\not=0$.

Considering the third case, we calculate the quotient $\frac{MP(\beta)}{MP(\beta^{-1})}$
 and check that this is not a Belyi function.

Considering the 4th case, we calculate the quotient $\frac{MP(\beta)}{MP(\beta^{-1})}$ and check
that it is a Belyi function. We calculate this quotient up to multiplication
$\beta\mapsto k\beta$, $k\in \kC$, $k\not=0$. The critical values are calculated for $r_3=s=1$
and are equal to $0$, $\infty$ and $\frac{525}{101838848}a(6040879+352815 a^2)$.

We have no more than $4$ conjugate solutions for our Belyi functions (we have 2 different dessins,
and 2 coordinate systems for each of them). Hence the 5th case can not give the desired solution.

\section{Working with the calculation results}
$j$-invariants of the found curves are irrational, hence we have the $2$ different curves,
so we have found at least two not isomorphic Belyi pairs with conditions stated in the beginning
(statement \ref{defdess}).
It completes the calculation. To simplify an answer, we represent it using new coordinates
(not the coordinates which used in calculations).

\begin{theorem}
\label{BPairs}
The Belyi pairs can be rewritten as follows:
Equation of the curve is $y^2=f(x)$. Belyi function looks as follows:
$$\beta=\frac{n_0-n_1+1}2 +
      y \sqrt{\frac{ (n_0-n_1)^2-2(n_0+n_1)+1}{4f}}.$$
(where if $\beta=P(x)+yQ(x)$ then $n_0=P^2-fQ^2$, $n_1=(P-1)^2-fQ^2$ or equivalently
$P=\frac{n_0-n_1+1}2$ and $Q^2=\frac{ (n_0-n_1)^2-2(n_0+n_1)+1}{4f}$).

\epsfbox{pict.1}\\
$\gamma=+45\sqrt{105}>0$ corresponds to the first dessin.

\epsfbox{pict.2}\\
$\gamma=-45\sqrt{105}<0$ corresponds to the second dessin
\end{theorem}

In terms of $\gamma=\pm 45\sqrt{105}$:
$$
n_0(x)=\frac{7^3 (39\gamma+17983)}{ 2^9 3^5 5^3}\cdot
\frac{(x+5)^3 (x-3)^5}{64x-105+\gamma}
$$
$$
n_1(x)= \frac{7^3 (39\gamma+17983)}{2^9 3^5 5^3} \cdot
\frac{(x^4-30x^2+40x+\frac{135\gamma-60825}{14})^2}{64x-105+\gamma}
$$
$$
f(x)=420 x^3 - (119+9\gamma) x^2+14(1515-\gamma)x+420(420-\gamma)
$$
\section{True shape.}
\begin{itemize}

\item
\epsfbox{pict.1}
If $j$-invariant is $j\approx 1315,640$ ($\gamma>0$) then we have a following true shape at the universal covering:
$$\epsfysize=5cm\epsfbox{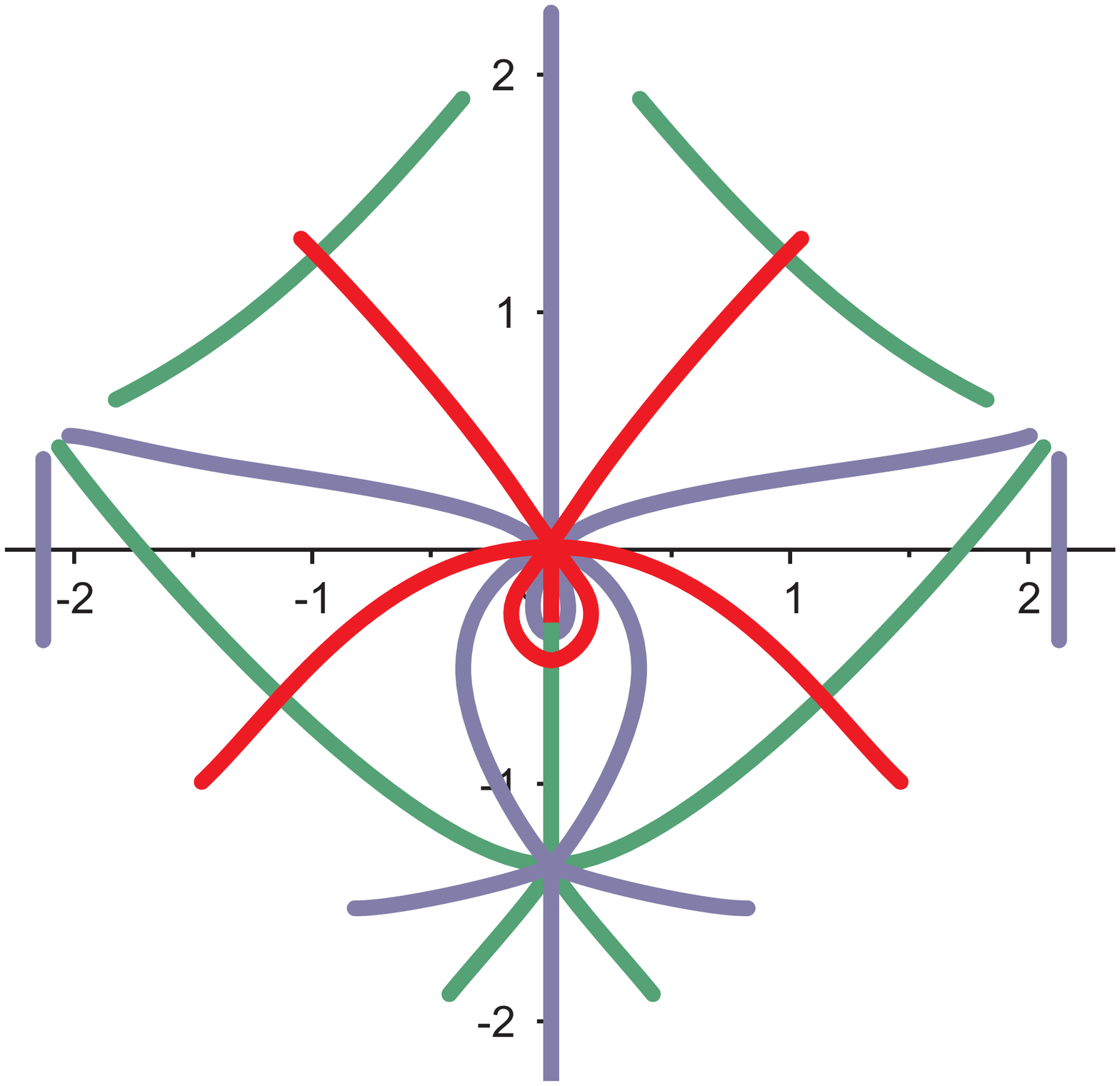}
\epsfysize=5cm\epsfbox{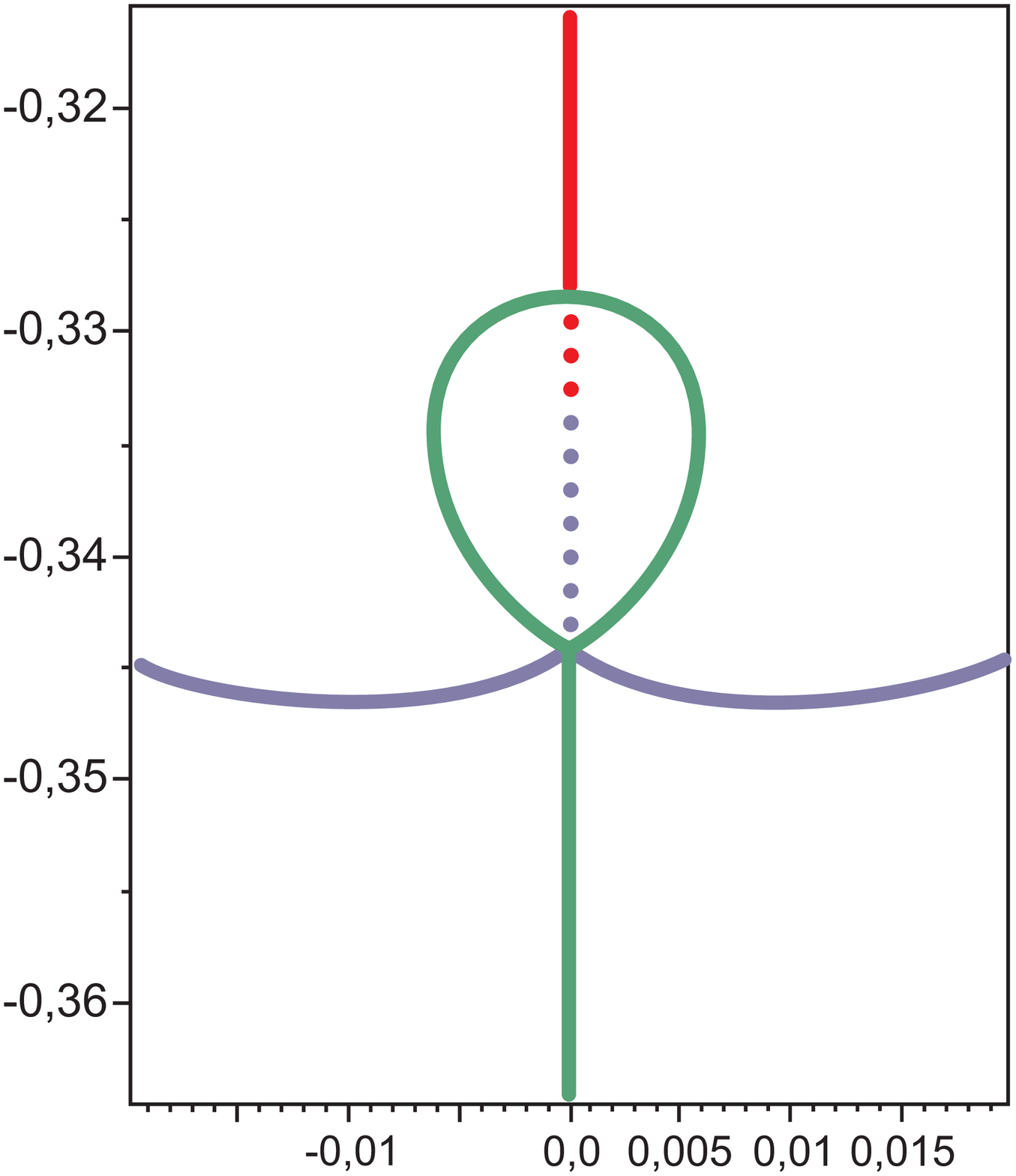}$$
\item
\epsfbox{pict.2}
If $j$-invariant is $j\approx 20,3167$ ($\gamma<0$) then we have a following true shape at the universal covering:
$$\epsfysize=5cm\epsfbox{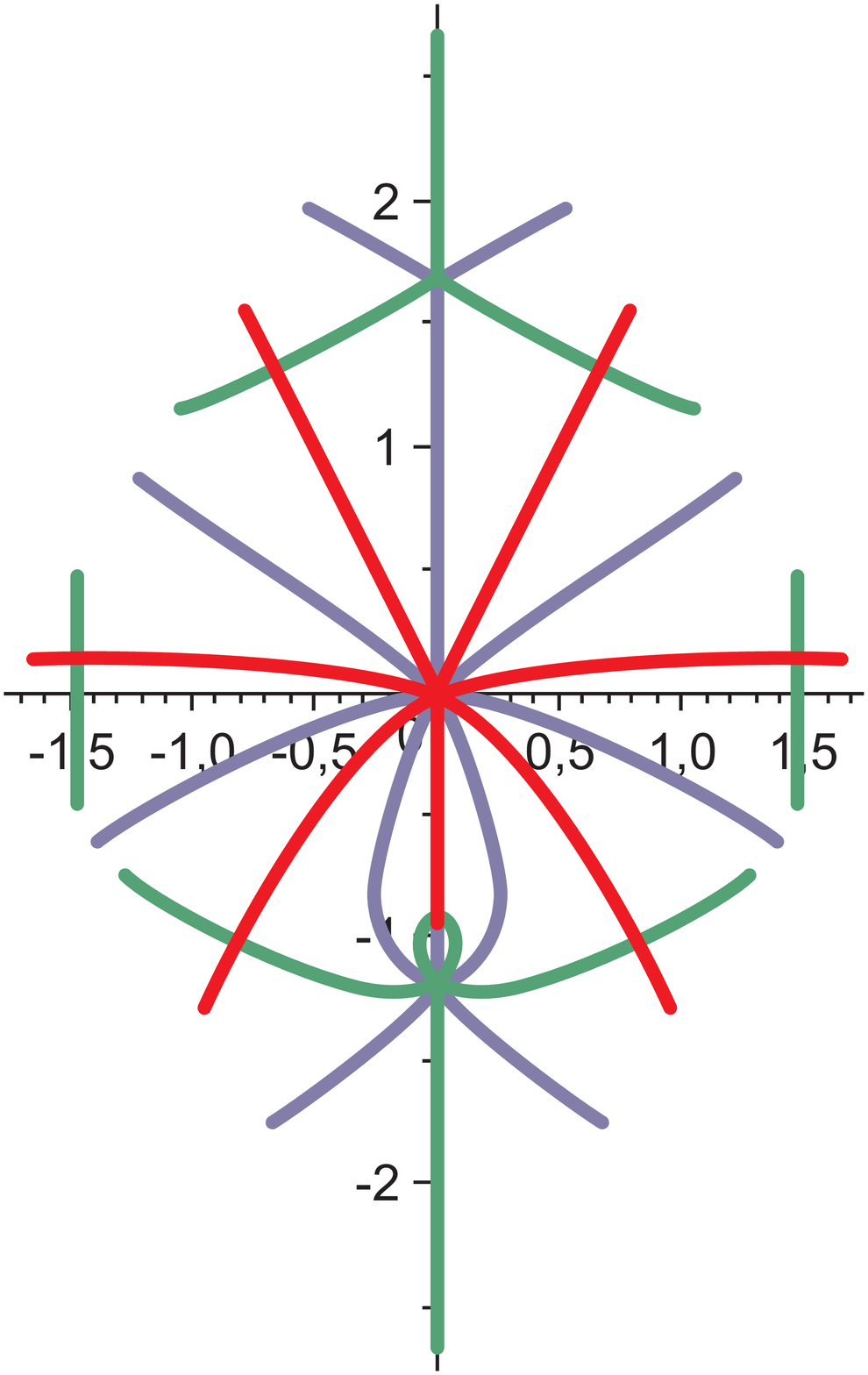}
\epsfysize=5cm\epsfbox{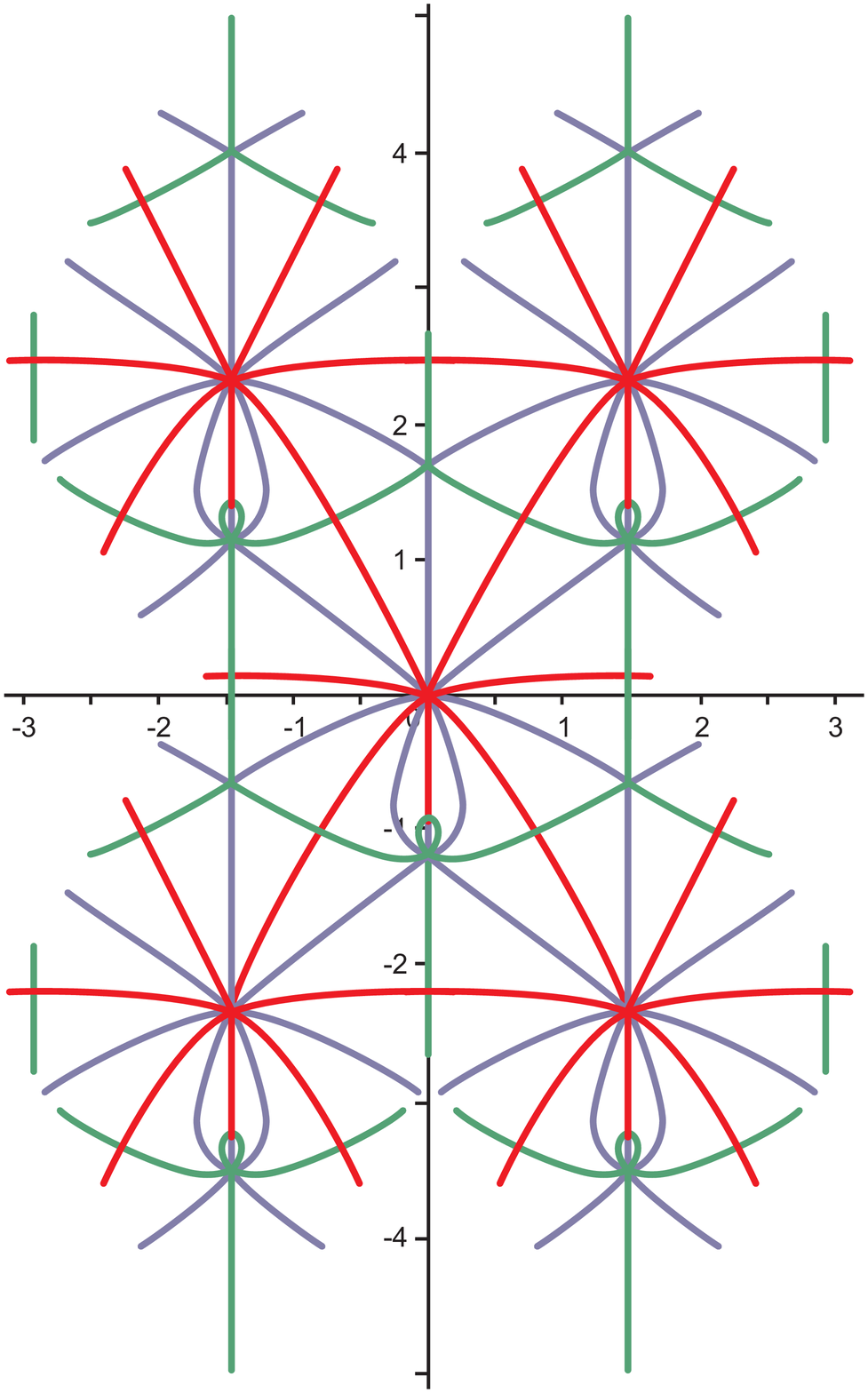}$$
\end{itemize}

\tableofcontents


\begin{thebibliography}{9}
\bibitem{LSchneps} Leila Schneps ''The Grothendieck Theory of Dessins D'enfants''
// Cambridge University Press, 1994,ISBN 0521478219, 9780521478212
\bibitem{Mul_Penk} M. Mulase and M. Penkava
''Ribbon Graphs, Quadratic Differentials on Riemann Surfaces, and
Algebraic Curves Defined over {$\overline{\mathbb Q}$}'' // Asian
Journal of Mathematics, Vol 2, number 4, 875-920 (1998).
\bibitem{Adr_Shab} N M Adrianov, G B Shabat,
''Belyi functions of dessins d'enfants of genus 2 with 4 edges''//
Russian Mathematical Surveys, 2005, 60:6, 1237-1239
\bibitem{LanZv} S.K.Lando, A.K.Zvonkin. Graphs on surfaces and their
applications. Encycl. of Math. Sciences, v.141, Springer, 2004.
\bibitem{maxplank} N.M. Adrianov, N.Ya. Amburg, V.A. Dremov, Yu.A.~Levitskaya,
E.M.~Kreines, Yu.Yu. Kochetkov, V.F.Nasretdinova, G.B.Shabat
''Catalog of dessins d'enfants with $\le 4$ edges''
//Max-Planck-Institut f\"ur Mathematik Preprint Series (MPI07-157), 2007
\end{thebibliography}
\end{document}